\RequirePackage{fix-cm} 
\documentclass[a4paper,twoside,12pt,reqno]{amsart}

\usepackage{fixltx2e}     

\usepackage[english]{babel}
\usepackage[latin1]{inputenc}

\usepackage{indentfirst,verbatim}

\usepackage{amsmath,amsfonts,amssymb,amsgen,amsbsy,eucal,mathrsfs,dsfont}
\usepackage{stmaryrd}

\usepackage[square,numbers]{natbib}
\usepackage{a4wide,verbatim}

\usepackage{epsfig,rotating,color}
\usepackage{multicol}
\usepackage{tikz,pgfplots}
\usepackage{graphicx}
\usepackage{mathtools}
\usepackage{textcomp}

\usepackage{caption}

\usepackage{hyperref}

\usepackage[square,numbers]{natbib}
\usepackage{systeme}

\hypersetup{
colorlinks=true, 
breaklinks=true, 
urlcolor= blue, 
linkcolor= blue, 
bookmarksopen=true, 
pdftitle={}, 
pdfauthor={}, 
pdfsubject={} 
}

\newcommand{\half}{{\textstyle{1\over2}}}
\newcommand{\third}{{\textstyle{1\over3}}}

\newcommand{\sixth}{{\textstyle{1\over6}}}
\newcommand{\fourth}{{\textstyle{1\over4}}}

\usepackage{amsthm}

\newtheorem{thm}{Theorem}
\newtheorem{lem}{Lemma}
\newtheorem{pro}{Proposition}

\theoremstyle{remark}
\newtheorem{rem}{Remark}

\usepackage{enumitem}
\newlist{steps}{enumerate}{1}
\setlist[steps, 1]{label = Step \arabic*:}

\newcommand{\eqdef}{\stackrel{\text{\tiny{def}}}{=}}

\title[Blow-up for the Serre--Green--Naghdi equations]{\bf On the blow-up scenario for some modified Serre--Green--Naghdi equations}

\author[GUELMAME]{Billel Guelmame}

\newcommand{\nfont}{\fontshape{n}\selectfont}

\address{({\nfont\textbf{Billel Guelmame}})  LJAD,  Inria,  Universit\'e C\^ote d'Azur, France.} 
\email{billel.guelmame@univ-cotedazur.fr}

\let\oldtocsection=\tocsection
 
\let\oldtocsubsection=\tocsubsection

\renewcommand{\tocsection}[2]{\hspace{0em}\oldtocsection {#1}{#2}}
\renewcommand{\tocsubsection}[2]{\hspace{2em}\oldtocsubsection{#1}{#2}}

\usepackage[pagewise]{lineno}
\begin{document}

\maketitle

\begin{abstract}
The present paper deals with a modified Serre--Green--Naghdi (mSGN) system that has been introduced by \citet{ClamondEtAl2017} to improve the dispersion relation. 
We present a precise blow-up scenario of the mSGN equations and we prove the existence of a class of solutions that develop singularities in finite time.
All the presented results hold, with the same proof, for the Serre--Green--Naghdi system with surface tension.
\end{abstract}

\medskip


\medskip


\tableofcontents

\section{Introduction}

\subsection{The equations}
Water waves are usually described by the Euler equations. Due to their complexity, other models have been proposed in various regimes.
For example, in the shallow water regime, it is assumed that the shallowness parameter $\sigma$ is small, where $\sigma$ is the ration of the mean water depth $\bar{h}$ to the wavelength $\iota$ (i.e., $\sigma = \bar{h}/\iota \ll 1$).
The Serre--Green--Naghdi equations are obtained by neglecting all the terms of order $\mathcal{O}(\sigma^4)$ from the water waves equations.

We consider in this paper the conservative modified Serre--Green--Naghdi system that was introduced by \citet{ClamondEtAl2017}
\begin{subequations}\label{SGN}
\begin{gather}
h_t\ +\,\left[\, h\,u\,\right]_x\ =\ 0, \label{SGNa} \\
\left[\,h\,u\,\right]_t\ +\,\left[\,h\,u^2\,+\, \half\, g\, h^2\, \label{SGN1b}
+\,\mathscr{R}\,\right]_x\ =\ 0,\\ \label{SGNc}
\mathscr{R}\, \eqdef\, \third\, \left( 1\, +\, {\textstyle \frac{3}{2}}\, \beta \right)\, h^3 \left(-u_{tx}\, -\, u\, u_{xx}\, +\, u_x^2 \right)\, 
-\, \half\, \beta\, g\, h^2\, \left( h\, h_{xx}\ +\ \half\, h_x^2 \right),
\end{gather}
\end{subequations}
where $h$ denotes the depth of the fluid, $u$ is the depth-averaged horizontal velocity, $g$ is the gravitational acceleration and $\beta$ is a free parameter. The classical Serre--Green--Naghdi system is recovered taking $\beta=0$.
The aim of this paper is to study the local (in time) well-posedness of \eqref{SGN}, to obtain a precise blow-up criterion and to build smooth solutions that develop singularities in finite time. 
 
Several modified Serre--Green--Naghdi equations have been derived and studied in the literature to optimise the linear dispersion,
some of them fail to conserve the energy and do not admit a variational principle, some other equations do not satisfy the Galilean invariance.
In order to improve the dispersion of the classical Serre--Green--Naghdi (cSGN) system conserving its desirable properties, \citet{ClamondEtAl2017} have modified the Lagrangian instead of modifying directly the cSGN equations to obtain the Lagrangian density (see Section \ref{sec:mod} for more details)
\begin{gather}\label{Lag}
\frac{\mathscr{L}_\beta}{\rho}\
=\ \half\, h\, u^2\ +\ \sixth\, \left( 1\, +\, {\textstyle \frac{3}{2}}\, \beta \right)\, h^3\, u_x^2\ -\ \half\, g\, h^2\ -\ \fourth\, \beta\, g\, h^2\,  h_x^2\ +\ \phi\, \left\{ h_t\ +\ [h\, u]_x \right\},
\end{gather}
where $\phi$ is a Lagrange multiplier.
The Euler--Lagrange equations of \eqref{Lag} lead to \eqref{SGN}.

Smooth solutions of the modified Serre--Green--Naghdi (mSGN) equations \eqref{SGN} satisfy the energy equation 
\begin{gather}\label{ene}
\mathscr{E}_t\ +\ \mathscr{Q}_x\ =\ 0,
\end{gather}
with 
\begin{gather}
\mathscr{E}\ \eqdef\ \half\, h\, u^2\ +\ \sixth\, \left( 1\, +\, {\textstyle \frac{3}{2}}\, \beta \right)\, h^3\, u_x^2\ +\ \half\, g\, \left(h\, -\, \bar{h}\right)^2\ +\ \fourth\, \beta\, g\, h^2\,  h_x^2, \\
\mathscr{Q}\ \eqdef\ u \left[ \mathscr{E}\ +\ \mathscr{R}\ +\ g\, h\, \left(h -\ \bar{h}\right) \right]\ +\ \half\, \beta\, g\, h^3\, h_x\, u_x,
\end{gather}
where $\bar{h}$ is the mean value of the depth $h$ of the fluid.
For $\beta >0$ and if $h$ is far from zero ($h \geqslant h_{min}>0$), the $H^1$ norms of both $u$ and $h-\bar{h}$ are controlled by the energy.
However, the energy of the cSGN equations ($\beta =0$) cannot control the $L^2$ norm of $h_x$.
This crucial property of mSGN \eqref{SGN} for $\beta>0$ is very important to build the small-energy solutions that develop singularities in finite time (Theorem \ref{thm:bu} below).
Several criteria have been proposed in \cite{ClamondEtAl2017} to chose the parameter $\beta$. 
For $\beta = 2/15 \approx 0.1333$, the dispersion relation of \eqref{SGN} coincide with the one of the full Euler system up order $4$ instead of order $2$ for the classical Serre--Green--Naghdi (Section \ref{sec:mod} below).
To optimise the decay of a particular solitary wave of \eqref{SGN} studied in \cite{ClamondEtAl2017}, one must take $\beta = {\textstyle \frac{2}{3}} \left( 12 \pi^{-2}-1 \right) \approx 0.1439$.
The value $\beta \approx 0.34560$ approximates the inner angle of the crest of the solitary wave solution to the exact angle of the limiting solitary wave ($120^{\circ}$).
In the present paper, we consider only the case $\beta>0$, and for the sake of simplicity we introduce
\begin{equation}
\alpha\ \eqdef\ 1\ +\ {\textstyle \frac{3}{2}}\, \beta.
\end{equation}

The mSGN equations on the form \eqref{SGN} contains some terms with high-order derivatives and a term with a time derivative in the definition of $\mathscr{R}$.
In order to obtain a simpler form of \eqref{SGN}, we introduce the  linear Sturm--Liouville operator
\begin{equation}\label{SL}
\mathcal{L}_h\ \eqdef\ h\ -\ \third\, \alpha\, \partial_x\, h^3\, \partial_x
\end{equation}
and we apply $\mathcal{L}_h^{-1}$ (the invertibility of $\mathcal{L}_h$ is proved in Lemma \ref{Inverseesitimates} below) on the equation \eqref{SGN1b} to obtain 
\begin{equation}\label{SGN1}
u_t\ +\ u\, u_x\ +  g\, h_x\ =\ - \mathcal{L}_h^{-1} \partial_x\, \left\{ {\textstyle \frac{2}{3}\, \alpha\, h^3\, u_x^2\ } +\ \third\, g\, h^3\, h_{xx}\ -\ \fourth\, \beta\, g\, h^2\, h_x^2 \right\}.
\end{equation}
In Lemma \ref{Inverseesitimates} below, we show that we gain one derivative with the operator $\mathcal{L}_h^{-1} \partial_x$, this is not enough to control the term $\third g h^3 h_{xx}$ in the right-hand side of \eqref{SGN1}.
To get rid of this term, we use \eqref{SL} to rewrite \eqref{SGN1} in the equivalent form
\begin{subequations}\label{SGN2}
\begin{align} \label{SGN2a}
h_t\ +\,\left[\, h\,u\,\right]_x\ &=\ 0, \\ \label{SGN2b}
u_t\ +\ u\, u_x\ + {\textstyle \frac{\alpha-1}{\alpha}}\, g\, h_x\ &=\ - \mathcal{L}_h^{-1}  \partial_x\, \left\{ {\textstyle \frac{2}{3}\, \alpha\, h^3\, u_x^2\ }  -\ \fourth\, \beta\, g\, h^2\, h_x^2\ +\ {\textstyle \frac{g}{2\, \alpha}}\, h^2 \right\}.
\end{align}
\end{subequations}
The left-hand side of \eqref{SGN2} is a symmetrisable $2 \times 2$ hyperbolic system and the right-hand side is a zero-order non-local term.
Then, the local well-posedness of \eqref{SGN2} in $H^s$ with $s \geqslant 2$ can be obtained following the proof of symmetrisable hyperbolic systems.
To the author's knowledge, the best blow-up criteria that have been obtained for those type of equations is ``if a singularity appears in finite time, then $\|h_x\|_{L^\infty} + \|u_x\|_{L^\infty}$ blows-up''.
This criteria does not show which term, or which slope blows-up. 
In this paper, we improve this criteria and we identify how exactly the blow-up occurs (Theorem \ref{thm:buc} and equation \eqref{pbc} below).
We also prove, using Riccati-type equations, the existence of a class of arbitrary small-energy initial data, such that the corresponding solutions develop singularities in finite time.

\subsection{Other similar equations}
In order to improve the dispersion relation of the classical Serre--Green--Naghdi equations, an interesting Whitham--Green--Naghdi (WGN) system have been proposed and studied in \cite{Duchene,Duchene2}.
The WGN system can be written as
\begin{subequations}\label{WGN}
\begin{gather}
h_t\ +\,\left[\, h\,u\,\right]_x\ =\ 0, \\
\left[ u - {\textstyle \frac{1}{3\, h}} \partial_x\, \mathrm{F} \left(h^3\, \partial_x\, \mathrm{F}\, u \right) \right]_t +\, g\, h_x\, +\, u\, u_x\, = \left[ {\textstyle \frac{u}{3\, h}} \partial_x\, \mathrm{F} \left(h^3\, \partial_x\, \mathrm{F}\, u \right) +\, \half\, h^2 \left( \partial_x\, \mathrm{F}\, u \right)^2  \right]_x,
\end{gather}
\end{subequations}
where $\mathrm{F}$ is the Fourier multiplier defined by
\begin{equation*}
\widehat{\mathrm{F} \varphi} (\xi)\ \eqdef\ {\textstyle \sqrt{\frac{3}{\bar{h}\, \xi\, \tanh \left(\bar{h}\, \xi \right)}\, -\, \frac{3}{\left(\bar{h}\, \xi \right)^2} }}\, \hat{\varphi}(\xi). 
\end{equation*}
The dispersion relation of the WGN system \eqref{WGN} is exactly the same as the one of the full Euler system. Smooth solutions of the WGN system conserve the energy $\mathcal{E}_t + \mathcal{D}_x=0$ with 
\begin{equation*}
\mathcal{E}\ \eqdef\ \half\, h\, u^2\ +\ \sixth\, h^3 \left( \partial_x\, \mathrm{F}\, u \right)^2\ +\ \half\, g\, \left(h\, -\, \bar{h}\right)^2.
\end{equation*}
Since the energy $\mathcal{E}$ does not control the $L^2$ norm of $h_x$, the results presented in this paper cannot be generalised directly for \eqref{WGN}. However, the WGN equations \eqref{WGN} deserve to be studied more in the future.

Other equations similar to \eqref{SGN} have been studied in the literature. For example, the Serre--Green--Naghdi equations with surface tension  \cite{dias,Lannes} 
\begin{subequations}\label{SGNsf}
\begin{gather}
h_t\ +\,\left[\, h\,u\,\right]_x\ =\ 0, \\
\left[\,h\,u\,\right]_t\ +\,\left[\,h\,u^2\,+\, \half\, g\, h^2\, 
+\,\mathscr{S}\,\right]_x\ =\ 0,\\ \label{SGNsfc}
\mathscr{S}\, \eqdef\, \third\, h^3 \left(-u_{tx}\, -\, u\, u_{xx}\, +\, u_x^2 \right)\, 
-\, \gamma \left( h\, h_{xx}\ -\ \half\, h_x^2 \right),
\end{gather}
\end{subequations}
where $\gamma>0$ is a constant (the surface tension coefficient divided by the density).
Instead of $\mathscr{S}$, the original equations involve 
\begin{equation*}
\tilde{\mathscr{S}}\, \eqdef\, \third\, h^3 \left(-u_{tx}\, -\, u\, u_{xx}\, +\, u_x^2 \right)\, 
-\ \gamma \left[ h \left(1 + h_x^2 \right)^{-3/2} h_{xx} +  \left(1 + h_x^2 \right)^{-1/2} \right].
\end{equation*}
Since we consider shallow water waves, any horizontal derivative is of order $\mathcal{O}(\sigma)$, using then Taylor expansion we obtain $\tilde{\mathscr{S}} = \mathscr{S} - \gamma + \mathcal{O}\left( \sigma^4 \right)$. Neglecting terms of order $\mathcal{O}\left(\sigma^4\right)$ we obtain \eqref{SGNsf}.
Smooth solutions of \eqref{SGNsf} conserve an $H^1$-equivalent energy. 
Weakly singular peakon travelling wave solutions of \eqref{SGNsf} have been studied in \cite{dut2,Mitsotakis} for the critical case $\gamma = g \bar{h}^2/3$.

Another similar system is the dispersionless regularised Saint-Venant (rSV) system proposed by \citet{ClamondDutykh2018a} that can be obtained replacing $\mathscr{R}$ in \eqref{SGN} by
\begin{equation*}
\mathscr{T}\, \eqdef\, \varepsilon\, h^3 \left(-u_{tx}\, -\, u\, u_{xx}\, +\, u_x^2 \right)\, 
-\, \varepsilon\, g\, h^2\, \left( h\, h_{xx}\ +\ \half\, h_x^2 \right),
\end{equation*}
where $\varepsilon >0$. The classical Saint-Venant equations are recovered by taking $\varepsilon=0$.
The rSV equations have been studied in \cite{ClamondDutykh2018a,liu2019well,PuEtAl2018} and have been generalised recently to regularise the Burgers equation \cite{rB} and the barotropic Euler equations \cite{guelmame2020Euler}. In \cite{liu2019well}, Liu et al. have proved the local well-posedness of the rSV equations and they derived smooth solutions that cannot exist globally in time. 
The proofs presented in this paper for the mSGN equations can be generalised for the SGN equations with weak surface tension \eqref{SGNsf}, for the rSV equations and also for the regularised barotropic Euler system \cite{guelmame2020Euler}.
The blow-up criterion proved here is more precise compared to the blow-up criterion in \cite{liu2019well}.
The key of the proof of the blow-up results in this paper is Lemma \ref{lem:PQ^2} below, a similar lemma have been proved in \cite{liu2019well} for a short time (shorter than the existence time in general). 
In this paper, we show, with a shorter (and different) proof, that the same result holds true as long as the smooth solution exists.


\subsection{Outline}
This paper is organised as follows. In Section \ref{sec:mod}, we present a brief derivation of the mSGN equations \eqref{SGN}. 
The main results are introduced in Section \ref{sec:mr}. 
In Section \ref{sec:pre}, we prove some useful lemmas.
Section \ref{sec:buc} is devoted to obtain the precise blow-up scenario of strong solutions of the mSGN equations.
In Section \ref{sec:bu},
we prove that some classical solutions cannot exist globally in time.

\section{Derivation}\label{sec:mod}

\subsection{Derivation}
In this section, we recall briefly the derivation of the model \eqref{SGN} presented in \cite{ClamondEtAl2017}. 
The classical Serre--Green--Naghdi equations can be derived from the Lagrangian density
\begin{gather}\label{Lag0}
\frac{\mathscr{L}_0}{\rho}\
=\ \half\, h\, u^2\ +\ \sixth\, h^3\, u_x^2\ -\ \half\, g\, h^2\  +\ \phi\, \left\{ h_t\ +\ [h\, u]_x \right\}.
\end{gather}
The Euler--Lagrange equations lead to 
\begin{equation*}
h_t\ +\,\left[\, h\,u\,\right]_x\ =\ 0, \qquad \qquad
u_t\ +\ u\, u_x\ +  g\, h_x\ =\ - {\textstyle \frac{1}{3\, h}} \left[ h^3 \left(-u_{tx}\, -\, u\, u_{xx}\, +\, u_x^2 \right) \right]_x.
\end{equation*}
Those equations describe long waves in the shallow water regime, thus, any horizontal and temporal derivative is of order $\mathcal{O}(\sigma)$.
This leads to 
\begin{equation*}
\left[u_t\ +\ u\, u_x\ +  g\, h_x\right]_x \, =\ \left[ - {\textstyle \frac{1}{3\, h}} \left[ h^3 \left(-u_{tx}\, -\, u\, u_{xx}\, +\, u_x^2 \right) \right]_x \right]_x\ =\ \mathcal{O}\left(\sigma^4\right).
\end{equation*}
In the Serre--Green--Naghdi model, all the terms of order $\mathcal{O}\left(\sigma^4\right)$ are neglected. 
\citet{ClamondEtAl2017} modified the Lagrangian \eqref{Lag0} by adding a neglected term to obtain the modified Lagrangian
\begin{equation*}
\frac{\tilde{\mathscr{L}}_\beta}{\rho}\ \eqdef\ \frac{\mathscr{L}_0}{\rho}\ +\ { \frac{\beta\, h^3}{12}} \left[u_t\ +\ u\, u_x\ +  g\, h_x\right]_x.
\end{equation*}
The choice of $h^3$ in the additional term is not necessary, one can replace $h^3$ by any increasing function of $h$ (see \cite{guelmame2020Euler} for more details). 
Using the conservation of the mass, one can write 
\begin{align*}
h^3 \left[u_t\ +\ u\, u_x\right]_x\ =\ [h^3\, u_x]_t\, +\ [h^3\, u_, u_x]_x\, +\ 3\, h^3\, u_x^2, \qquad 
h^3\, h_{xx}\ =\ [h^3\, h_x]_x\, -\ 3\, h^2\, h_x^2.
\end{align*}
This leads to $\tilde{\mathscr{L}}_\beta = \mathscr{L}_\beta + [\cdots]_t + [\cdots]_x$, where $\mathscr{L}_\beta$ is the Lagrangian density defined in \eqref{Lag}.
This means that $\mathscr{L}_\beta \equiv \tilde{\mathscr{L}}_\beta = \mathscr{L}_0 + \mathcal{O}\left(\sigma^4\right)$.
Deriving the Lagrangian density $\mathscr{L}_\beta$ we obtain \eqref{SGN}.
Thus, the mSGN system \eqref{SGN} is a suitable modification of the cSGN system that admits good properties and conserves a better energy compared to the cSGN equations for $\beta>0$.

\subsection{The dispersion relation}
The choice of the free parameter $\beta$ depends on the desired property. For example, in order to improve the dispersion relation of \eqref{SGN}, one can linearise \eqref{SGN} around the constant state $\left(\bar{h},0\right)$ and consider traveling waves on the form $\cos \left\{ k(x - c t) \right\}$.
The dispersion relation then becomes
\begin{equation*}
\frac{c^2}{g\, \bar{h}}\ =\ \frac{2\, +\, \beta \left(k\, \bar{h}\right)^2}{2\, + \left( {\textstyle \frac{2}{3}}\, +\,  \beta\right) \left(k\, \bar{h}\right)^2}\ =\ 1\ -\ \third \left(k\, \bar{h}\right)^2\ +\ \left( {\textstyle \frac{1}{9} +\, \frac{\beta}{6}} \right) \left(k\, \bar{h}\right)^4\ +\ \cdots.
\end{equation*}
This must be compared with the exact relation
\begin{equation*}
\frac{c^2}{g\, \bar{h}}\ =\ \frac{\tanh \left(k\, \bar{h}\right)}{k\, \bar{h}}\ =\ 1\ -\ \third \left(k\, \bar{h}\right)^2\ +\ {\textstyle \frac{2}{15}}\left(k\, \bar{h}\right)^4\ +\ \cdots .
\end{equation*}
Hence, in order to to have the same dispersion relation up to the order $k^4$, one must take $\beta = 2/15$.

Other criteria exits for choosing the parameter $\beta$ (see \cite{ClamondEtAl2017} for more details).

\section{Main results}\label{sec:mr}
The first result of this paper is the local well-posedness of the system \eqref{SGN2} in the Sobolev space 
\begin{equation}
H^s(\mathds{R})\ \eqdef\ \left\{f,\, \|f\|_{H^s(\mathds{R})}^2\ \eqdef\ \int_{\mathds{R}} \left( 1\, +\, |\xi|^2 \right)^s\, |\hat{f}(\xi)|^2\, \mathrm{d}\xi\ <\ +\infty \right\}
\end{equation}
where $s \geqslant 2$ is a real number.
\begin{thm}\label{thm:ex}
Let $\beta > 0$, $\bar{h}>0$ and $s \geqslant 2$, then, for any $(h_0-\bar{h},u_0) \in H^s(\mathds{R})$ satisfying $\inf_{x\, \in\, \mathds{R}}h_0(x)>0$ there exists $T>0$ and $(h-\bar{h},u) \in C^0([0,T],H^s(\mathds{R})) \cap C^1([0,T],H^{s-1}(\mathds{R}))$ a unique solution of \eqref{SGN2} such that 
\begin{equation}\label{ndc}
\inf_{(t,x)\, \in\, [0,T] \times \mathds{R}}\ h(t,x)\ >\ 0.
\end{equation}
Moreover, the solution satisfies the conservation of the energy
\begin{equation}\label{eneC}
\frac{\mathrm{d}}{\mathrm{d}t}\, \int_\mathds{R} \left( \half\, h\, u^2\ +\ \sixth\, \alpha\, h^3\, u_x^2\ +\ \half\, g\, \left(h\, -\, \bar{h}\right)^2\ +\ \fourth\, \beta\, g\, h^2\,  h_x^2 \right)\, \mathrm{d}x\ =\ 0.
\end{equation}
\end{thm}
\begin{rem}
The solution given in Theorem \ref{thm:ex} depends continuously on the initial data, i.e., If $(h_0-\bar{h},u_0), (\tilde{h}_0-\bar{h},\tilde{u}_0) \in H^s$, such that $h_0,\tilde{h}_0 \geqslant h^* >0 $, and $t \leqslant \min \{T, \tilde{T} \}$, then there exists a constant $C(\|(\tilde{h}-\bar{h},\tilde{u})\|_{L^\infty([0,t],H^{s})},\|(h-\bar{h},u)\|_{L^\infty([0,t],H^{s})})>0$, such that 
\begin{equation}
\|(h-\tilde{h},u-\tilde{u})\|_{L^\infty([0,t],H^{s-1})}\ \leqslant\ C\, \|(h_0-\tilde{h}_0,u_0-\tilde{u}_0)\|_{H^s}.
\end{equation}
\end{rem}
The proof of Theorem \ref{thm:ex} is classic and omitted in this paper. See \cite{guelmame2020Euler,Israwi2011,Lannes,liu2019well} and Theorem 3 of \cite{guelmame2020rSV} for more details.
\begin{rem}
If $(h-\bar{h},u) \in C^0([0,T],H^s(\mathds{R}))$ and $h$ satisfies \eqref{ndc} for some $T>0$, then Theorem \ref{thm:ex} ensures that the solution can be extended over $[0,T]$. In other words, if $T_{max}$ is the maximum time existence of the solution, then, we have the blow-up criterion 
\begin{equation}\label{buc0}
T_{max}\ <\ +\infty\ \implies\ \liminf_{t\to T_{max}}\, \inf_{x\, \in\, \mathds{R}}\, h(t,x)\ =\ 0\quad \mathrm{or}\quad  \limsup_{t\to T_{max}}\,  \|(h-\bar{h},u)\|_{H^s}\ =\ +\infty . 
\end{equation}
\end{rem}

\noindent
The blow-up criterion \eqref{buc0} can be improved as in \cite{guelmame2020Euler,guelmame2020rSV,liu2019well} to 
\begin{equation}\label{buc0.5}
T_{max}\ <\ +\infty\ \implies\ \liminf_{t\to T_{max}}\, \inf_{x\, \in\, \mathds{R}}\, h(t,x)\ =\ 0\quad \mathrm{or}\quad  \limsup_{t\to T_{max}}\,  \|(h_x,u_x)\|_{L^\infty}\ =\ +\infty . 
\end{equation}
The last blow-up criterion ensures that if $h>0$ is far from zero, then, the blow-up will appear on the $L^\infty$ norm of $u_x$ or $h_x$. 
This result is improved in this paper, and we claim the two more precise criteria for blow-up mechanism.
\begin{thm}\label{thm:buc}
Let $ \beta >0$, and let $T_{max}$ be the maximum time existence of the solution given by Theorem \ref{thm:ex}, then
\begin{equation}\label{buc1}
T_{max}\, <\, +\infty \ \implies \ \liminf_{t\to T_{max}}\, \inf_{x\, \in\, \mathds{R}}\, h(t,x)\ =\ 0\quad \mathrm{or}\quad
\left\{\begin{aligned} \liminf_{t\to T_{max}}\, \inf_{x\, \in\, \mathds{R}}\, &u_x(t,x)\ =\ -\infty, \\ \mathrm{and}& \\ \limsup_{t\to T_{max}}\, \|h_x&(t,x)\|_{L^\infty}\, =\, +\infty, \end{aligned}\sysdelim. . \right.
\end{equation}
which is equivalent to second criterion
\begin{equation}\label{buc2}{\small
T_{max}\, <\, +\infty \ \implies \ 
\limsup_{t\to T_{max}}\, \|u_x(t,x)\|_{L^\infty}\, =\, +\infty\ \mathrm{and}\  
\left\{\begin{aligned} \liminf_{t\to T_{max}}\, \inf_{x\, \in\, \mathds{R}}\, &h(t,x)\ =\ 0, \\ \mathrm{or}& \\ \limsup_{t\to T_{max}}\, \|h_x(t,&x)\|_{L^\infty}\, =\, +\infty. \end{aligned}\sysdelim. . \right.}
\end{equation}
\end{thm}

Since $H^1 \hookrightarrow L^\infty$ and the energy \eqref{eneC} is conserved, then, $|h-\bar{h}|$ is controlled by the energy of the initial data (see Proposition \ref{pro:ene} below). This ensures that if the initial energy is small enough then $h$ is uniformly far from zero ($\min_{t,x} h >0$).
In this case, the blow-up criterion \eqref{buc1} becomes 
\begin{equation}\label{pbc}
T_{max}\, <\, +\infty \ \implies \ 
\liminf_{t\to T_{max}}\, \inf_{x\, \in\, \mathds{R}}\,  u_x(t,x)\ =\ -\infty
\quad
\mathrm{and}
\quad
\limsup_{t\to T_{max}}\,  \|h_x\|_{L^\infty}\ =\ +\infty.
\end{equation} 
This blow-up criterion ensures that if a blow-up occurs, then $u_x$ goes to $-\infty$. Since $h$ is far from zero, then the conservation of the mass implies that the material derivative of the free surface goes to $+\infty$.
However, it is not clear if $h_x$ blows-up on $- \infty$ or $+ \infty$.  
The following theorem shows that both scenarios are possible.
\begin{thm}\label{thm:bu}
For any $T>0$ and $K \in \left]0,{\textstyle \frac{g \sqrt{\beta}}{3 \sqrt{2}}} \bar{h}^3 \right[$, there exist
\begin{itemize}
\item  $(h_0-\bar{h},u_0) \in C^\infty_c(\mathds{R})$ satisfying $\int_\mathds{R} \mathscr{E}_0\, \mathrm{d}x \leqslant K$ such that the corresponding solution of \eqref{SGN2} blows-up in finite time $T_{max} \leqslant T$ and 
\begin{equation*}
\inf_{[0,T_{max}[ \times \mathds{R}}  u_x(t,x)\ =\ -\infty,
\quad
\sup_{[0,T_{max}[ \times \mathds{R}}  h_x(t,x)\ =\ +\infty,
\quad
\inf_{[0,T_{max}[ \times \mathds{R}}  h_x(t,x)\ >\ -\infty.
\end{equation*}
\item  $(\tilde{h}_0-\bar{h},\tilde{u}_0) \in C^\infty_c(\mathds{R})$ satisfying $\int_\mathds{R} \tilde{\mathscr{E}}_0\, \mathrm{d}x \leqslant K$ such that the corresponding solution of \eqref{SGN2} blows-up in finite time $\tilde{T}_{max} \leqslant T$ and 
\begin{equation*}
\inf_{[0,\tilde{T}_{max}[ \times \mathds{R}}  \tilde{u}_x(t,x)\ =\ -\infty,
\quad
\inf_{[0,\tilde{T}_{max}[ \times \mathds{R}}  \tilde{h}_x(t,x)\ =\ -\infty,
\quad
\sup_{[0,\tilde{T}_{max}[ \times \mathds{R}}  \tilde{h}_x(t,x)\ <\ +\infty.
\end{equation*}
\end{itemize}
\end{thm}

\begin{rem}
All the proofs presented in this paper work also for the SGN equations with surface tension \eqref{SGNsf}, for the rSV equations \cite{ClamondDutykh2018a} and also for the regularised barotropic Euler system \cite{guelmame2020Euler}.
\end{rem}

\section{Preliminaries}\label{sec:pre}
In this section, we recall some classical estimates and we prove some lemmas that are needed to prove Theorem \ref{thm:buc} and Theorem \ref{thm:bu}.

\begin{lem}(\cite{constantin2002initial})
Let $F \in {C}^{\tilde{m}+2}(\mathds{R})$ with $F(0)=0$ and $0 \leqslant s \leqslant \tilde{m}$, then there exists a continuous function $\tilde{F}$, such that for all $f \in H^s \cap W^{1,\infty}$ we have
\begin{equation}\label{Composition2}
\|F(f)\|_{H^s}\ \leqslant\ \tilde{F} \left( \|f\|_{W^{1,\infty}} \right)\, \|f\|_{H^s}.
\end{equation}
\end{lem}
Let $\Lambda$ be defined such that $\widehat{\Lambda f}=(1+\xi^2)^\frac{1}{2} \hat{f}$, then we have the following estimate.
\begin{lem}(\cite{kato1988commutator})
Let $[A,B] \eqdef AB-BA$ be the commutator of the operators $A$ and $B$. If $r\, \geqslant\, 0$, then $\exists C>0$ such that
\begin{align}
\|f\, g\|_{H^r}\ &\leqslant\ C\, \left( \|f\|_{L^\infty}\, \|g\|_{H^r}\ +\ \|f\|_{H^r}\, \|g\|_{L^\infty}\right), \label{Algebra} \\
\left\| \left[ \Lambda^r,\, f \right]\, g  \right\|_{L^2}\ &\leqslant\ C\, \left( \|f_x\|_{L^\infty}\, \|g\|_{H^{r-1}}\ +\ \|f\|_{H^r}\, \|g\|_{L^\infty} \right). \label{Commutator}
\end{align}
\end{lem}

Now, we recall the invertibility of the operator $\mathcal{L}_h$ defined in \eqref{SL} and that we gain two derivatives with $\mathcal{L}_h^{-1}$.
\begin{lem}\label{Inverseesitimates}
Let $\alpha>0$ and $0<  h \in W^{1,\infty}$ with $h^{-1} \in L^\infty$, then the operator $\mathcal{L}_h$ is an isomorphism from $H^2$ to $L^2$ and   $\exists C_1 = C_1\left(\alpha,s,  \left\| h^{-1} \right\|_{L^\infty}, \|h - \bar{h}\|_{W^{1,\infty}} \right)>0$, 
$C_2=C_2\left(\alpha,\left\| h^{-1} \right\|_{L^\infty}, \|h\|_{L^{\infty}} \right)>0$ such that 
\begin{enumerate}
\item  If $ s \geqslant 0$, then
\begin{subequations}\label{estimate1}
\begin{align}
\left\|\mathcal{L}_h^{-1}\, \partial_x\, \psi \right\|_{H^{s+1}}\  
&\leqslant\ C_1\, \left( \left\|\psi \right\|_{H^s}\ +\ \left\|h\, -\, \bar{h} \right\|_{H^s}\, 
\left\|\mathcal{L}_h^{-1}\, \partial_x\, \psi \right\|_{W^{1,\infty}}\right),\\
\left\|\mathcal{L}_h^{-1}\, \phi \right\|_{H^{s+1}}\  
&\leqslant\ C_1\, \left( \left\|\phi \right\|_{H^s}\ +\ \left\|h\, -\, \bar{h} \right\|_{H^s}\, 
\left\|\mathcal{L}_h^{-1}\, \phi \right\|_{W^{1,\infty}} \right).
\end{align}
\end{subequations}
\item   If $ s \geqslant 0$, then  
\begin{subequations}\label{estimate2}
\begin{align}
\left\|\mathcal{L}_h^{-1}\, \partial_x\, \psi \right\|_{H^{s+1}}\  
&\leqslant\ C_1\, \left\|\psi \right\|_{H^s}\, \left(1 +\ \left\|h\, -\, \bar{h} \right\|_{H^s} \right), \\
\left\|\mathcal{L}_h^{-1}\, \phi \right\|_{H^{s+1}}\  
&\leqslant\ C_1\,  \left\|\phi \right\|_{H^s}\, \left(1 +\ \left\|h\, -\, \bar{h} \right\|_{H^s} \right).
\end{align}
\end{subequations}
\item If $\phi \in {C}_{\mathrm{lim}} \eqdef \{f\in {C}(\mathds{R}),\ f(+\infty), f(-\infty) \in \mathds{R} \}$, then $\mathcal{L}_h^{-1}\, \phi$ is well defined and 
\begin{equation}\label{estimate3}
\left\|\mathcal{L}_h^{-1}\, \phi\right\|_{W^{2,\infty}}\ \leqslant\ C_2\, \left\|\phi\right\|_{L^\infty}. 
\end{equation}
\item If $\psi \in {C}_{\mathrm{lim}} \cap L^1$, then 
\begin{equation}\label{estimate4}
\left\|\mathcal{L}_h^{-1}\, \partial_x \psi\right\|_{W^{1,\infty}}\ \leqslant\ C_2\, \left( \left\|\psi\right\|_{L^\infty}\ +\ \left\|\psi\right\|_{L^1} \right). 
\end{equation}
\item If $h_x \in L^2$, then
\begin{subequations}\label{estimate56}
\begin{align}\label{estimate5}
\left\| \mathcal{L}_{h}^{-1}\, \partial_x\, \psi \right\|_{H^1}\, +\, \left\| \mathcal{L}_{h}^{-1}\, \psi \right\|_{H^1}\, &\leqslant\ C_2 \, \|\psi\|_{L^2},\\ \label{estimate6}
\left\| \mathcal{L}_{h}^{-1}\,  \psi \right\|_{W^{1,\infty}}\, \leqslant\
\left\| \mathcal{L}_{h}^{-1}\,  \psi \right\|_{H^2}\,  &\leqslant\ C_2 \left[1\, +\, \left\|h_x \right\|_{L^2}^2 \right]\left\| \psi \right\|_{L^2}.
\end{align}
\end{subequations}
\end{enumerate}
\end{lem}

\proof 
The proof of the invertibility of $\mathcal{L}_h$ and the estimates \eqref{estimate1},  \eqref{estimate2}, \eqref{estimate3}, \eqref{estimate4} and \eqref{estimate5} can be found in \cite{liu2019well}. It remains only to prove \eqref{estimate6}.
Using the definition of $\mathcal{L}_h$ we obtain
\begin{align*}
\partial_x^2\,\mathcal{L}_h^{-1}\, \psi\ 
=&\ \partial_x\, h^{-3} \mathcal{L}_h^{-1} \left[h^4\, \partial_x\, \mathcal{L}_h^{-1} \psi\, -\, {\textstyle \frac{\alpha}{3} \partial_x\, h^3\, \partial_x\, h^3\, \partial_x\, \mathcal{L}_h^{-1} \psi } \right]\\
=&\ \partial_x\, h^{-3} \mathcal{L}_h^{-1}\, \partial_x \left[h^4\,  \mathcal{L}_h^{-1} \psi\, -\, {\textstyle \frac{\alpha}{3} h^3\, \partial_x\, h^3\, \partial_x\, \mathcal{L}_h^{-1} \psi } \right] -\, 4\, \partial_x\, h^{-3} \mathcal{L}_h^{-1}\, h^3\, h_x \mathcal{L}_h^{-1} \psi\\
=&\ \partial_x\, h^{-3} \mathcal{L}_h^{-1}\, \partial_x h^3\, \psi\,  -\, 4\, \partial_x\, h^{-3} \mathcal{L}_h^{-1}\, h^3\, h_x \mathcal{L}_h^{-1} \psi\\
=&\  -3\, h^{-2}\, h_x \left[\mathcal{L}_{h}^{-1}\, \partial_x\, h^3\, \psi  - 4\,
\mathcal{L}_{h}^{-1}\, \left[ h^3\, h_x\, \mathcal{L}_h^{-1}\, \psi\right] \right]\
+\ h^{-3}\, \partial_x\, \mathcal{L}_{h}^{-1}\, \partial_x\, h^3\, \psi\\
&-\ 4\,
h^{-3}\, \partial_x\, \mathcal{L}_{h}^{-1}\, \left[\, h^3\, h_x\, \mathcal{L}_h^{-1}\, \psi\right].
\end{align*}
Using \eqref{estimate5} and the embedding $H^1 \hookrightarrow  L^\infty$ we obtain 
\begin{align*}
\left\| \partial_x^2\, \mathcal{L}_{h}^{-1}\,  \psi \right\|_{L^2}\,  
\leqslant &\  C_2 \left\|h_x \right\|_{L^2} \left[ \left\|\mathcal{L}_{h}^{-1}\, \partial_x\, h^3\, \psi \right\|_{H^1}\ +\ \left\| \mathcal{L}_{h}^{-1}\, \left[ h^3\, h_x\, \mathcal{L}_h^{-1}\, \psi\right] \right\|_{H^1} \right]\\
&+\ C_2 \left\| \partial_x\, \mathcal{L}_{h}^{-1}\, \partial_x\, h^3\, \psi\ \right\|_{L^2}\ +\ C_2 \left\|\partial_x\, \mathcal{L}_{h}^{-1}\, \left[ h^3\, h_x\, \mathcal{L}_h^{-1} \psi \right]  \right\|_{L^2}\\
\leqslant &\ C_2 \left\|h_x \right\|_{L^2} \left[ \left\| \psi \right\|_{L^2}\ +\ \left\|h_x \right\|_{L^2} \left\|  \mathcal{L}_h^{-1}\, \psi\right\|_{H^1} \right]\\
&+\ C_2 \left\| \psi \right\|_{L^2}\ +\ C_2 \left\|h_x \right\|_{L^2} \left\|  \mathcal{L}_h^{-1}\, \psi\right\|_{H^1}\\
\leqslant &\ C_2 \left[1\, +\, \left\|h_x \right\|_{L^2}^2 \right]\left\| \psi \right\|_{L^2}.
\end{align*}
This with \eqref{estimate5} imply \eqref{estimate6}.
\qed

The $\mathscr{R}$ defined in \eqref{SGNc} contains some terms with two order derivatives, using \eqref{SGN2}, we can write $\mathscr{R}$ without those high order derivatives involving the operator $\mathcal{L}_h^{-1}$. Then, using the previous lemma we show that the norm $\|\mathscr{R}\|_{L^\infty}$ is controlled by  $\|(h,u,h^{-1})\|_{L^\infty}$.
\begin{lem}\label{lem:inqR}
Let $(h-\bar{h},u)$ be a smooth solution of \eqref{SGN2}, then for any $T < T_{max}$, there exists $C=C(\beta,\bar{h},\|(h,u,h^{-1})\|_{L^\infty([0,T]\times \mathds{R})})>0$, such that
\begin{equation}\label{inqR}
\left\|\mathscr{R}\right\|_{L^\infty([0,T] \times \mathds{R})}\ \leqslant\ C.
\end{equation}
\end{lem}

\proof
From the definition of $\mathcal{L}_h$, we obtain that
\begin{equation}\label{Psi}
\left( 1\ +\ \third\, \alpha\, h^3\, \partial_x\, \mathcal{L}_h^{-1}\, \partial_x \right)\, \Psi\ =\ h^3\, \partial_x\, \mathcal{L}_h^{-1}\, \left( h \int_{-\infty}^x h^{-3}\, \Psi \right)
\end{equation}
for any smooth function $\Psi$, such that $\Psi(\pm \infty)=0$. Using \eqref{SGNc}, \eqref{SGN2b} and \eqref{Psi} we obtain
\begin{align}
\mathscr{R}\ &=\ - \third\, \alpha\, h^3\, \partial_x\, \left[u_t\ +\ u\, u_x\ +\ {\textstyle \frac{\alpha-1}{\alpha} }\, g\, h_x \right]\ +\ {\textstyle \frac{2}{3}}\, \alpha\, h^3\, u_x^2\ -\ {\textstyle \frac{1}{4}}\, \beta\, g\, h^2 h_x^2 \nonumber\\ 
&=\ \left( 1\ +\ \third\, \alpha\, h^3\, \partial_x\, \mathcal{L}_h^{-1}\, \partial_x \right) \left\{ {\textstyle \frac{2}{3}\, \alpha\, h^3\, u_x^2\ }  -\ \fourth\, \beta\, g\, h^2\, h_x^2\ +\ g\, {\textstyle \frac{h^2\, -\, \bar{h}^2}{2\,\alpha}} \right\}\ -\ g\, {\textstyle \frac{h^2\, -\, \bar{h}^2}{2\,\alpha}} \label{RR}\\
&=\ \left( 1\ +\ \third\, \alpha\, h^3\, \partial_x\, \mathcal{L}_h^{-1}\, \partial_x \right) \left\{ {\textstyle \frac{2}{3}\, \alpha\, h^3\, u_x^2\ }  -\ \fourth\, \beta\, g\, h^2\, h_x^2 \right\}\, +\ {\textstyle \frac{g}{3}}\, h^3\, \partial_x\, \mathcal{L}_h^{-1}\, \left\{ h\, h_x \right\}  \nonumber\\
&=\ h^3\, \partial_x\, \mathcal{L}_h^{-1}\, \left( h \int_{-\infty}^x \left( {\textstyle \frac{2}{3}\, \alpha\, u_x^2\ }  -\ \fourth\, \beta\, g\, h^{-1}\, h_x^2 \right) \right)\ +\ {\textstyle \frac{g}{3}}\, h^3\, \partial_x\, \mathcal{L}_h^{-1}\, \left\{ h\, h_x \right\}
\end{align}
Using the conservation of the energy \eqref{eneC} we obtain 
\begin{align*}
\left\| \int_{-\infty}^x\! \left( {\textstyle \frac{2}{3}\, \alpha\, u_x^2\, }  -\, \fourth\, \beta\, g\, h^{-1}\, h_x^2 \right) \right\|_{L^\infty}\!
&\leqslant\ \left\| {\textstyle \frac{2}{3}\, \alpha\, u_x^2\, }  -\, \fourth\, \beta\, g\, h^{-1}\, h_x^2 \right\|_{L^1}\
\leqslant\ C_3.
\end{align*}
Then, the inequality \eqref{inqR} follows directly from \eqref{estimate3} and \eqref{estimate6}.
\qed

Since we are considering the Serre--Green--Naghdi equations on the form \eqref{SGN2} instead of \eqref{SGN1}, it is more convenient to use the following Riemann invariants\footnote{Those quantities are constants along the characteristics if the right-hand side of \eqref{SGN2} is zero.}
 $R$, $S$ and their corresponding speeds of characteristics $\lambda$, $\mu$ 
\begin{align}
R\ &\eqdef\ u\ +\ 2\, \sqrt{{\textstyle \frac{\alpha-1}{\alpha}\, g\, h}}, \qquad 
\lambda\ \eqdef\ u\ +\ \sqrt{{\textstyle \frac{\alpha-1}{\alpha}\, g\, h}}, \\
S\ &\eqdef\ u\ -\ 2\, \sqrt{{\textstyle \frac{\alpha-1}{\alpha}\, g\, h}}, \qquad 
\mu\ \eqdef\ u\ -\ \sqrt{{\textstyle \frac{\alpha-1}{\alpha}\, g\, h}},
\end{align}
rather that the Riemann invariants of the classical Saint-Venant system.
Then, the system \eqref{SGN2} can be rewritten as 
\begin{subequations}\label{SGNR}
\begin{align}
R_t\ +\ \lambda\, R_x\ &=\ - \mathcal{L}_h^{-1}  \partial_x\, \left\{ {\textstyle \frac{2}{3}\, \alpha\, h^3\, u_x^2\ }  -\ \fourth\, \beta\, g\, h^2\, h_x^2\ +\ {\textstyle \frac{g}{2\, \alpha}}\, h^2 \right\}, \\
S_t\ +\ \mu\, S_x\ &=\ - \mathcal{L}_h^{-1}  \partial_x\, \left\{ {\textstyle \frac{2}{3}\, \alpha\, h^3\, u_x^2\ }  -\ \fourth\, \beta\, g\, h^2\, h_x^2\ +\ {\textstyle \frac{g}{2\, \alpha}}\, h^2 \right\}.
\end{align}
\end{subequations}
Defining
\begin{align*}
P\ \eqdef\ R_x\ =\ u_x\ +\ \sqrt{{\textstyle \frac{\alpha-1}{\alpha}\, g}}\, h^{-\half} h_x,\\
Q\ \eqdef\ S_x\ =\ u_x\ -\ \sqrt{{\textstyle \frac{\alpha-1}{\alpha}\, g}}\, h^{-\half} h_x,
\end{align*}
we have 
\begin{equation}\label{uh_x}
u_x\ =\ \frac{P\ +\ Q}{2}, \qquad h_x\ =\  \frac{\sqrt{\alpha}\, h^{\frac{1}{2}}}{2\, \sqrt{(\alpha -1)\, g}} \left( P\ -\ Q \right).
\end{equation}
Let the characteristics $X_{a},Y_{a}$ starting from $a$ defined as the solutions of the ordinary differential equations
\begin{align}
\frac{\mathrm{d}}{\mathrm{d}t}\, X_{a}(t)\ &=\ \lambda(t,X_{a}(t)), \qquad X_{a}(0)\ =\ a\\
\frac{\mathrm{d}}{\mathrm{d}t}\, Y_{a}(t)\ &=\ \, \mu(t,Y_{a}(t)), \qquad\ Y_{a}(0)\ =\ a
\end{align}
Differentiation \eqref{SGNR} with respect to $x$, and using \eqref{RR} we obtain the Ricatti-type equations
\begin{subequations}\label{Ricatti}
\begin{align}\label{Ricatti1}
\frac{\mathrm{d}^\lambda}{\mathrm{d}t}\, P\ \eqdef\ P_t\ +\ \lambda\, P_x\ =\ - {\textstyle \frac{3}{8}}\, P^2\ +\ {\textstyle \frac{3}{8}}\, Q^2\ +\ P\, Q\  -\ 3\, \alpha^{-1}\, h^{-3}\, \mathscr{R},\\ \label{Ricatti2}
\frac{\mathrm{d}^\mu}{\mathrm{d}t}\, Q\ \eqdef\ Q_t\ +\ \mu\, Q_x\ =\ - {\textstyle \frac{3}{8}}\, Q^2\ +\ {\textstyle \frac{3}{8}}\, P^2\ +\ P\, Q\  -\ 3\,\alpha^{-1}\, h^{-3}\, \mathscr{R},
\end{align}
\end{subequations}
where $\frac{\mathrm{d}^\lambda}{\mathrm{d}t}, \frac{\mathrm{d}^\mu}{\mathrm{d}t}$ denote the derivatives along the characteristics with the speeds $\lambda$ and $\mu$ respectively.

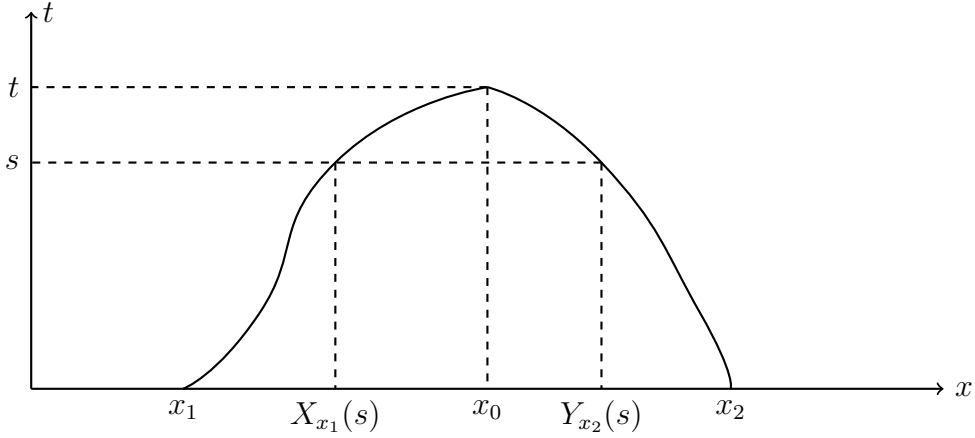
\begin{figure}[!ht]
\begin{tikzpicture}[thick, transform canvas={scale=1}, shift={(0,-5)}]
\draw[->] (-6,0) -- (6,0) node[right]{$x$};
 \draw[->] (-6,0) -- (-6,5) node[right]{$t$};

\draw [black] plot [smooth, tension=1] coordinates { (3.2,0) (2.8,1) (1.5,3) (0,4)};

\draw [black] plot [smooth, tension=1] coordinates { (-4,0) (-3,1) (-2,3) (0,4)};

\draw[dashed] (0,4) -- (-6,4) node[left]{$t$};
\draw[dashed] (0,4) -- (0,0) node[below]{$x_0$};
\draw[dashed] (1.5,3) -- (-6,3) node[left]{$s$};
\draw[dashed] (1.5,3) -- (1.5,0) node[below]{$Y_{x_2}(s)$};
\draw[dashed] (-2,3) -- (-2,0) node[below]{$X_{x_1}(s)$};


\draw[] (3.2,0) -- (3.2,0) node[below]{$x_2$};
\draw[] (-4,0) -- (-4,0) node[below]{$x_1$};

\end{tikzpicture}
\vspace{5.7cm}
\caption{Characteristics.}
\label{XY}
\end{figure}

A key point to prove Theorem \ref{thm:buc} and Theorem \ref{thm:bu} is to control the term $P^2$ in the Ricatti equation \eqref{Ricatti2} and the term $Q^2$ in \eqref{Ricatti1}.
For that purpose, we prove in the following lemma that the integral of $P^2$ along the $X$ characteristics and the integral of $Q^2$ along the $Y$ characteristics are bounded.
\begin{lem}\label{lem:PQ^2}
Let $\beta > 0$, $\bar{h}>0$ and $(h_0-\bar{h},u_0) \in H^2$ initial data satisfying $\inf_{x\, \in\, \mathds{R}}h_0(x)>0$ and let $(h-\bar{h},u)$ be the corresponding solution of \eqref{SGN2} given by Theorem \ref{thm:ex}, let also $t\in [0,T_{max}[$, then, there exist {\small
$$ A\left(\beta,\bar{h},\|(h,u,h^{-1})\|_{L^\infty([0,t]\times \mathds{R})}, \int\! \mathscr{E}\, \mathrm{d}x \right)>0,\ 
B\left(\beta,\bar{h},\|(h,u,h^{-1})\|_{L^\infty([0,t]\times \mathds{R})}, \int\! \mathscr{E}\, \mathrm{d}x\right)>0,$$ }
\noindent
such that for any $x_2 \in \mathds{R}$,  and for $x_1 \in ]-\infty, x_2[$ the solution of $X_{x_1}(t)=Y_{x_2}(t)$ (see Figure \ref{XY}) we have 
\begin{align}\label{PQ^2}
\int_0^t Q(s,X_{x_1}(s))^2\, \mathrm{d}s\ +\ \int_0^t P(s,Y_{x_2}(s))^2\, \mathrm{d}s\ \leqslant\ A\, t\ +\ B.
\end{align}
\end{lem}

\begin{rem}
A similar result have been proved for the so-called variational wave equation with $A=0$ and $B$ depends only on the energy of the initial data \cite{GlasseyEtAl1996}.
For the mSGN, additional terms appear, and a uniform (on time) bound cannot be obtained for large data.
\end{rem}

\proof 
Defining 
\begin{align*}
\mathscr{B}_1\ &\eqdef\  
\sqrt{\textstyle \frac{\alpha-1}{\alpha}\, g\, h}\, \left(\half\, h\, u^2\ +\ \half\, g\, \left(h\, -\, \bar{h} \right)^2 \right)\
-\  u\, \left( \mathscr{R}\ +\ g\, h\, \left(h -\ \bar{h}\right) \right), \\
\mathscr{B}_2\ &\eqdef\  
\sqrt{\textstyle \frac{\alpha-1}{\alpha}\, g\, h}\, \left(\half\, h\, u^2\ +\ \half\, g\, \left(h\, -\, \bar{h} \right)^2 \right)\
+\  u\, \left( \mathscr{R}\ +\ g\, h\, \left(h -\ \bar{h}\right) \right), 
\end{align*}
and using \eqref{inqR}, one can prove that the quantity $\|\mathscr{B}_1\|_{L^\infty} + \|\mathscr{B}_2\|_{L^\infty}$ is bounded. One can easily check that 
\begin{equation}\label{ED}
\lambda\, \mathscr{E}\ -\ \mathscr{D}\ =\ 
{\textstyle \frac{\sqrt{6\, \alpha\, \beta\, g\, h^7}}{12}}\, Q^2\ +\ \mathscr{B}_1, \qquad -\mu\, \mathscr{E}\ +\ \mathscr{D}\ =\ {\textstyle \frac{\sqrt{6\, \alpha\, \beta\, g\, h^7}}{12}}\, P^2\ +\ \mathscr{B}_2 .
\end{equation}
Since 
$$ \lambda\ -\ \mu\ =\ 2\, \sqrt{{\textstyle \frac{\alpha-1}{\alpha}\, g\, h}}\
\geqslant\ 2\, \sqrt{{\textstyle \frac{\alpha-1}{\alpha}\, g}}\, \|h^{-1}\|_ {L^\infty}^{-\frac{1}{2}}\ >\ 0, $$
then $x_1 < x_2$.
Integrating \eqref{ene} on the set $\left\{ (s,x),\ s \in [0,t], X_{x_1}(s) \leqslant x \leqslant Y_{x_2}(s)  \right\}$, using the divergence theorem, the energy equation \eqref{eneC} and \eqref{ED} one obtains \eqref{PQ^2}. \qed

\section{Blow-up criteria}\label{sec:buc}

The aim of this section is to prove Theorem \ref{thm:buc}, for that purpose, we consider $s \geqslant 2$ and $(h-\bar{h},u)$ the solution of \eqref{SGN2} given by Theorem \ref{thm:ex} with the initial data $(h_0-\bar{h},u_0)$ and we start by the following lemmas.
\begin{lem}\label{lem:equiv}
For $T_{max}< +\infty$, we consider the following properties 
\begin{subequations}
\begin{gather} \label{i}
\sup_{(t,x)\, \in\, [0,T_{max}[ \times \mathds{R}} u_x(t,x)\ <\ +\infty,\\ \label{ii}
\inf_{(t,x)\, \in\, [0,T_{max}[ \times \mathds{R}}  h(t,x)\ >\ 0,\\ \label{iii}
\left\| (h,u) \right\|_{L^\infty([0,T_{max}[ \times \mathds{R})}\ <\ +\infty.
\end{gather}
\end{subequations}
Then, \eqref{i} $\iff$  \eqref{ii} and \eqref{ii} $\implies$ \eqref{iii}.
\end{lem}

\proof
The proof of \eqref{ii} $\implies$ \eqref{iii} follows directly from the conservation of the energy \eqref{eneC} and the embedding $H^1 \hookrightarrow  L^\infty$.

Let the characteristic $Z_{a}$ starting from $a$ defined as the solutions of the ordinary differential equation
\begin{align}
\frac{\mathrm{d}}{\mathrm{d}t}\, Z_{a}(t)\ &=\ u(t,Z_{a}(t)), \qquad Z_{a}(0)\ =\ a.
\end{align}
Denoting the derivatives along the characteristics with speed $u$ by $\frac{\mathrm{d}^u}{\mathrm{d}t}$ and using the conservation of the mass \eqref{SGN2a} one obtains
\begin{equation}
\frac{\mathrm{d}^u}{\mathrm{d}t}\, h\ \eqdef\ h_t\ +\ u\, h_x\  =\ - u_x\, h, \qquad \implies \qquad h\ \geqslant\ \left( \inf_{x\, \in\, \mathds{R}} h_0 \right)\, \mathrm{e}^{- \sup_{t,x}\,  u_x(t,x)\, T_{max}}.
\end{equation}
The proof of \eqref{i} $\implies$  \eqref{ii} follows directly from the last inequality.

It only remains to prove the converse (\eqref{ii} $\implies$  \eqref{i}). Using the Young inequality 
\begin{equation}\label{Young}
\pm a\, b\ \leqslant\ {\textstyle \frac{3}{8}}\, a^2\ +\ {\textstyle \frac{2}{3}}\, b^2, 
\end{equation}
integrating \eqref{Ricatti1}, \eqref{Ricatti2} along the characteristics, and using \eqref{inqR}, \eqref{PQ^2} one obtains the existence of $\tilde{A}>0, \tilde{B}>0$ which depend on $\beta,\bar{h}$, $\|(h,u,h^{-1})\|_{L^\infty([0,T]\times \mathds{R})}$ and $\int\! \mathscr{E}\, \mathrm{d}x$, such that
\begin{subequations}\label{PQ<}
\begin{align}\label{P<}
P(t,X_{x_1}(t))\ \leqslant\ P_0(x_1)\ +\ \tilde{A}\, t\ +\ \tilde{B} \qquad \forall (t,x_1) \in [0,T]\times \mathds{R},\\ \label{Q<}
Q(t,Y_{x_2}(t))\ \leqslant\ Q_0(x_2)\ +\ \tilde{A}\, t\ +\ \tilde{B} \qquad \forall (t,x_2) \in [0,T]\times \mathds{R}.
\end{align}
\end{subequations}
The last inequalities imply that 
\begin{equation}\label{iiii}
\sup_{(t,x)\, \in\, [0,T_{max}[ \times \mathds{R}} P(t,x)\ <\ +\infty, \qquad \mathrm{and} \qquad \sup_{(t,x)\, \in\, [0,T_{max}[ \times \mathds{R}} Q(t,x)\ <\ +\infty,
\end{equation}
then, \eqref{i} follows directly from \eqref{uh_x}.
\qed

\begin{lem}\label{lem:equiv2}
For $T_{max}< +\infty$, we consider the following properties 
\begin{subequations}
\begin{gather} \label{1}
\inf_{(t,x)\, \in\, [0,T_{max}[ \times \mathds{R}} u_x(t,x)\ >\ -\infty,\\ \label{2}
\left\| h_x \right\|_{L^\infty([0,T_{max}[ \times \mathds{R})}\ <\ +\infty.
\end{gather}
\end{subequations}
Then, \eqref{ii} $\implies$ $\left(\eqref{1} \iff \eqref{2} \right)$.
\end{lem}

\proof 
We suppose that \eqref{ii} and \eqref{1} are satisfied. Using Lemma \ref{lem:equiv} one obtains that $\|u_x\|_{L^\infty}$ is bounded. Then \eqref{2} follows directly from \eqref{iiii} and
\begin{equation}
h_x\ =\  
\frac{\sqrt{\alpha}\, h^{\frac{1}{2}}}{ \sqrt{(\alpha -1)\, g}} \left( u_x\ -\ Q \right)\ =\
\frac{\sqrt{\alpha}\, h^{\frac{1}{2}}}{ \sqrt{(\alpha -1)\, g}} \left( P\ -\ u_x \right).
\end{equation}
To prove the converse, we suppose that \eqref{ii} and \eqref{2} are satisfied, then, using the Young inequality $\pm ab \geqslant -\frac{1}{2} a^2 - \frac{1}{2} b^2$, \eqref{Ricatti1} and 
\begin{equation}
P^2\ =\ Q^2\ +\ 4\, \sqrt{\textstyle \frac{\alpha-1}{\alpha}\, g}\, h^{-\frac{1}{2}} u_x\, h_x\ =\ Q^2\ +\ 2\, \sqrt{\textstyle \frac{\alpha-1}{\alpha}\, g}\, h^{-\frac{1}{2}}\, h_x\, (P\, +\, Q)
\end{equation}
one obtains
\begin{align*}
\frac{\mathrm{d}^\lambda}{\mathrm{d}t}\, P\ 
&\geqslant\ - {\textstyle \frac{7}{8}}\, P^2\ -\ {\textstyle \frac{1}{8}}\, Q^2\ -\ 3\, \alpha^{-1}\, h^{-3}\, \mathscr{R}\\
&=\ -Q^2\ 
-\ {\textstyle \frac{7}{4}}\, \sqrt{\textstyle \frac{\alpha-1}{\alpha}\, g}\, h^{-\frac{1}{2}}\, h_x\, P\ 
-\ {\textstyle \frac{7}{4}}\, \sqrt{\textstyle \frac{\alpha-1}{\alpha}\, g}\, h^{-\frac{1}{2}}\, h_x\, Q\
-\ 3\, \alpha^{-1}\, h^{-3}\, \mathscr{R}\\
&\geqslant\ -\ {\textstyle \frac{7}{4}}\, \sqrt{\textstyle \frac{\alpha-1}{\alpha}\, g}\, h^{-\frac{1}{2}}\, h_x\, P\ -\ {\textstyle \frac{3}{2}}\, Q^2\
-\ {\textstyle \frac{49}{32}\, \frac{\alpha-1}{\alpha}\, } g\, h^{-1}\, h_x^2\
-\ 3\, \alpha^{-1}\, h^{-3}\, \mathscr{R}.
\end{align*}
Using \eqref{2}, \eqref{inqR}, \eqref{PQ^2} and Gronwall lemma, we obtain that $\inf_{t,x} P > - \infty$. Using again \eqref{2} we obtain \eqref{1}.
\qed

\vspace{0.5cm}

Now, we can prove Theorem \ref{thm:buc}. Note that Lemma \ref{lem:equiv} implies the equivalence between \eqref{buc1} and \eqref{buc2}.
Then, it only remains to prove \eqref{buc1}.
Step 1 is devoted to prove the blow-up criterion \eqref{buc0.5}.
The proof of \eqref{buc1} is given in Step 2.

\textit{Proof of Theorem \ref{thm:buc}.}

\textbf{\textit{Step 1:}} In order to prove \eqref{buc0.5}, we suppose that $\left\| (h_x,u_x) \right\|_{L^\infty([0,T_{max}[ \times \mathds{R})} < +\infty$, \eqref{ii} and we prove that if  $T_{max}<+\infty$, then $\left\| (h-\bar{h},u) \right\|_{L^\infty([0,T_{max}[, H^s(\mathds{R})} < +\infty$ which contradicts with the definition of $T_{max}$. 
For that purpose, we define
\begin{gather*}
W\ \eqdef\ (h\, -\, \bar{h},u)^\top \qquad
A(W)\ \eqdef\ \begin{pmatrix} \frac{\alpha-1}{\alpha}\, g & 0 \\ 0 & h \end{pmatrix}, \qquad   
B(W)\ \eqdef\ \begin{pmatrix} u & h \\ \frac{\alpha-1}{\alpha}\, g & u \end{pmatrix},
\\
\mathscr{F}(W)\ \eqdef\ \begin{pmatrix} 0 \\ - \mathcal{L}_h^{-1}  \partial_x\, \left\{ {\textstyle \frac{2}{3}\, \alpha\, h^3\, u_x^2\ }  -\ \fourth\, \beta\, g\, h^2\, h_x^2\ +\ {\textstyle \frac{g}{2\, \alpha}}\, h^2 \right\} \end{pmatrix},
\end{gather*}
the system \eqref{SGN2} becomes
\begin{equation}\label{system}
W_t\ +\ B(W)\, W_x\ =\ \mathscr{F}(W).
\end{equation}
Let $(\cdot,\cdot)$ be the scalar product in $L^2$ and 
$ E(W)\ \eqdef\ \left(\Lambda^s\, W,\ A\, \Lambda^s\, W \right)$.
Since $AB$ is a symmetric matrix, straightforward calculations with \eqref{system} show that
\begin{align}\label{energyterms}
E(W)_t\ 
=&\ -2\, \left( [\Lambda^s, B]\, W_x,\, A\, \Lambda^s W \right)\ -\ 2\, \left(B\, \Lambda^s  W_x,\, A\, \Lambda^s W \right) \nonumber \\
&\ - 2\, \left( \Lambda^s \mathscr{F},\, A\, \Lambda^s W \right)\ +\ \left( \Lambda^s W,\, A_t\, \Lambda^s W \right) \nonumber \\
=&\ -2\, \left( [\Lambda^s, B]\, W_x,\, A\, \Lambda^s W \right)\ +\ \left( \Lambda^s  W,\, (A\, B)_x\, \Lambda^s W \right) \nonumber \\
&\ - 2\, \left( \Lambda^s \mathscr{F},\, A\, \Lambda^s W \right)\ +\ \left( \Lambda^s W,\, A_t\, \Lambda^s W \right) 
\end{align}
Defining $\bar{B} \eqdef B(\bar{h},0)$, and using \eqref{Commutator}, \eqref{Composition2} one obtains
\begin{align*} 
\left| \left( [\Lambda^s, B]\, W_x, A\, \Lambda^s W \right)\ \right|\ 
&\leqslant\ C\, \|A\|_{L^\infty} \|W\|_{H^s} \left( \|B_x\|_{L^\infty} \|W_x\|_{H^{s-1}}\, +\, \|B-\bar{B}\|_{H^s} \|W_x\|_{L^\infty}\!  \right) \\
&\leqslant\ \tilde{C}_1\, \|W\|_{H^s}^2,
\end{align*}
where $\tilde{C}_1$ is a positive constant that depends on $\|W\|_{W^{1,\infty}}$ and $\|h^{-1}\|_{L^\infty}$. Using the conservation of the mass \eqref{SGN2a} one obtains that 
\begin{equation}
\left| \left( \Lambda^s  W,\, (A\, B)_x\, \Lambda^s W \right)\right|\ +\ \left|   \left( \Lambda^s W,\, A_t\, \Lambda^s W \right)  \right|\
\leqslant\  \tilde{C}_2\, \|W\|_{H^s}^2.
\end{equation}
From the energy conservation \eqref{eneC}, it is clear that 
\begin{equation}
\left\| {\textstyle \frac{2}{3}\, \alpha\, h^3\, u_x^2\ }  -\ \fourth\, \beta\, g\, h^2\, h_x^2\right\|_{L^1}\, +\, \left\| g\, {\textstyle \frac{h^2-\bar{h}^2}{2\, \alpha}} \right\|_{L^2}\ \leqslant\ \tilde{C}_3\, 
\end{equation}
Then, using \eqref{estimate1}, \eqref{estimate4}, \eqref{estimate5}, \eqref{Algebra} and \eqref{Composition2} we obtain
\begin{equation}
\|\mathscr{F}\|_{H^s}\ \leqslant\ \tilde{C}_4\, \|W\|_{H^s}.
\end{equation}
Since $\|(h,h^{-1})\|_{L^\infty}$ is bounded, then, $\|W\|_{H^s} \leqslant \tilde{C}_5 E(W)$. Combining all the estimates above, one obtains 
\begin{equation}
E(W)_t\ \leqslant\ \tilde{C}\, E(W),
\end{equation}
where $\tilde{C}$ does not depend on $\|W\|_{H^s}$. Then, Gronwall lemma implies that 
$$\left\| (h-\bar{h},u) \right\|_{L^\infty([0,T_{max}[, H^s(\mathds{R})} \leqslant \tilde{C}_5 \|E(W)\|_{L^\infty ([0,T_{max}[)} < +\infty.$$
This ends the proof of \eqref{buc0.5}.

\textbf{\textit{Step 2:}} It remains to prove \eqref{buc1}. We suppose that $T_{max}<+\infty$ and \eqref{ii} is satisfied. The blow-up criterion \eqref{buc0.5} (previous step) insures that 
\begin{equation*}
\limsup_{t\to T_{max}}\,  \|u_x\|_{L^\infty}\  =\ +\infty 
\qquad
\mathrm{or}
\qquad
\limsup_{t\to T_{max}}\,  \|h_x\|_{L^\infty}\ =\ +\infty .
\end{equation*}
Lemma \ref{lem:equiv} implies that $u_x$ is bounded from above, then 
\begin{equation*}
\liminf_{t\to T_{max}}\, \inf_{x\, \in\, \mathds{R}}\,  u_x(t,x)\ =\ -\infty
\qquad
\mathrm{or}
\qquad
\limsup_{t\to T_{max}}\,  \|h_x\|_{L^\infty}\ =\ +\infty .
\end{equation*}
Finally, Lemma \ref{lem:equiv2}, insures that if one of the quantities above blows-up, then the other one should also blow-up at the same time. Then 
\[
\pushQED{\qed} 
\liminf_{t\to T_{max}}\, \inf_{x\, \in\, \mathds{R}}\,  u_x(t,x)\ =\ -\infty
\qquad
\mathrm{and}
\qquad
\limsup_{t\to T_{max}}\,  \|h_x\|_{L^\infty}\ =\ +\infty .\qedhere
\popQED
\]

\section{Blow-up results}\label{sec:bu}
The goal of this section is to prove Theorem \ref{thm:bu}, for that purpose, we consider smooth solutions with small energy. 
In the following proposition we prove that if the energy is small enough, then, the quantity $\|(h,u,h^{-1})\|_{L^\infty}$ is uniformly bounded.
\begin{pro}\label{pro:ene}
For $\beta > 0$, $\bar{h}>0$, let $E$ be a positive number such that
\begin{align*}
0\ <\ E\ <\ {\textstyle \frac{g\, \sqrt{\beta}}{3\, \sqrt{2}}}\, \bar{h}^3, 
\end{align*}
Defining 
\begin{gather*}
h_{min} \eqdef \bar{h}\, -\, \sqrt{\textstyle \frac{3\, \sqrt{2}\, E}{g\, \sqrt{\beta}\, \bar{h}}}, \qquad 
h_{max} \eqdef \bar{h}\, +\, \sqrt{\textstyle \frac{3\, \sqrt{2}\, E}{g\, \sqrt{\beta}\, \bar{h}}}, \\ 
u_{max} \eqdef -u_{min} \eqdef (3/\alpha)^\frac{1}{4}\, \sqrt{E}/h_{min}.
\end{gather*}
Then, for any $(h-\bar{h},u) \in H^1$, if $\int \mathscr{E}\, \mathrm{d}x \leqslant E$, we have
\begin{gather}\label{ubound}
0\ <\ h_{min}\ \leqslant\ h\ \leqslant h_{max}\ < 2\, \bar{h}, \qquad u_{min}\ \leqslant\ u\ \leqslant\ u_{max}, 
\end{gather}
\end{pro}

\begin{rem}\label{rem:constants}
The conservation of the energy \eqref{eneC} and Proposition \ref{pro:ene} insure that if the initial data satisfy $\int \mathscr{E}_0\, \mathrm{d}x \leqslant E$, then, as long the the solution exist, the quantity $\|(h,u,h^{-1})\|_{L^\infty}$ is bounded by a constant that depends only on $g, \gamma, \bar{h}$ and $E$.
Thence, all the constants given in \eqref{estimate3}, \eqref{estimate4}, \eqref{inqR}, \eqref{PQ^2} and \eqref{PQ<} are universal and do not depend on the initial data and the solution.
\end{rem}

\vspace{0.15cm}
\textit{Proof of Proposition \ref{pro:ene}.}
The Young inequality $\frac{1}{2} a^2 + \frac{1}{2} b^2 \geqslant \pm ab$ implies that
\begin{align*} 
E\ \geqslant\ \int_\mathds{R} \mathscr{E}\, \mathrm{d}y\ 
\geqslant &\ 
\int_\mathds{R} \left( \half\, g\, \left(h\, -\, \bar{h}\right)^2\ +\ \fourth\, \beta\, g\, h^2\,  h_x^2 \right)\, \mathrm{d}x\\
\geqslant &\ g\, \sqrt{\textstyle \frac{\beta}{2}}\, \left( \int_{-\infty}^x (h\, -\, \bar{h})\, h\, h_x\, \mathrm{d}y\ -\ \int^{+\infty}_x (h\, -\, \bar{h})\, h\, h_x\, \mathrm{d}y  \right)\\
\geqslant &\ {\textstyle \frac{g}{3}}\, \sqrt{\textstyle \frac{\beta}{2}}\, \left( h\ -\ \bar{h} \right)^2 \left( 2\, h\ +\ \bar{h} \right)\\
\geqslant &\ {\textstyle \frac{g}{3}}\, \sqrt{\textstyle \frac{\beta}{2}}\, \bar{h}\, \left( h\ -\ \bar{h} \right)^2, 
\end{align*}
which implies that $h_{min} \leqslant h \leqslant h_{max}$. Doing the same estimates with $u$ one obtains
\begin{align*} 
E\ \geqslant\ \int_\mathds{R} \mathscr{E}\, \mathrm{d}y\ 
\geqslant &\  \int_\mathds{R} \left( \half h\, u^2\ +\ \sixth\, \alpha\, h^3\, u_x^2 \right)\, \mathrm{d}y\\
\geqslant &\ {\textstyle \sqrt{\frac{\alpha}{3}}\, h_{min}^2}\, \left(  \int_{-\infty}^x u\, u_x\, \mathrm{d}y\ -\ \int^{+\infty}_x u\, u_x\, \mathrm{d}y \right)\\
\geqslant &\ {\textstyle \sqrt{\frac{\alpha}{3}}\, h_{min}^2}\, |u|^2,
\end{align*}
the last inequality ends the proof of $u_{min} \leqslant u \leqslant u_{max}$. \qed

\vspace{0.3cm}
\textit{Proof of Theorem \ref{thm:bu}.}
Since the proofs of the two parts of Theorem \ref{thm:bu} are the same, we only prove the first part.

Let $T>0$, and let $\tilde{A},\tilde{B}$ be the constants given in \eqref{PQ<}.
From \eqref{inqR} we obtain that $|\alpha^{-1}\, h^{-3} \mathscr{R}| \leqslant \tilde{C}^2$ for some $\tilde{C}>0$.
If the initial data satisfy $\int \mathscr{E}_0\, \mathrm{d}x \leqslant E$, then, the constants $\tilde{A},\tilde{B}$ and $\tilde{C}$ are universal and depend only on $g, \beta, \bar{h}$ and $E$ (Remark \ref{rem:constants}).
We choose $\tilde{T}$ and $\tilde{D}$ such that 
\begin{equation}
0\ <\ \tilde{T}\ \leqslant\ T, \qquad 3\, \tilde{C}\, \tilde{T}\ \leqslant\ \frac{\pi}{4}, \qquad  \tilde{D}\ \eqdef\ 4\, \max \left\{\tilde{A}^2,\tilde{B}^2,\tilde{C}^2 \right\},
\end{equation}
and we choose the initial data $(h_0-\bar{h},u_0)\in C^\infty_c(\mathds{R})$ such that there exist $x_1 \in \mathds{R}$ satisfying
\begin{equation}
\int_{\mathds{R}} \mathscr{E}_0\, \mathrm{d}x\ \leqslant\ E,  \quad Q_0\ \equiv\ 0, \quad P_0(x_1)\ <\ -\ 2\, \sqrt{ \tilde{D}}\, \left( \tilde{T}\, +\, 1 \right), \quad P_0(x_1)\ <\ - \frac{8}{\tilde{T}}.
\end{equation}
Let $t< \min \{ \tilde{T},T_{max} \}$, then \eqref{Ricatti2} with the Young inequality $PQ \geqslant -{\textstyle \frac{3}{8}} P^2 -{\textstyle \frac{2}{3}} Q^2$ imply that
\begin{equation}
\frac{\mathrm{d}^\mu}{\mathrm{d}t}\, Q\ \geqslant\ - {\textstyle \frac{25}{24}}\, Q^2\  -\ 3\, \alpha^{-1}\, h^{-3}\, \mathscr{R}\ \geqslant\ - 3\, \left( Q^2\  +\ \tilde{C}^2 \right).
\end{equation}
The last inequality with \eqref{Q<} imply that for all $x_1 \in \mathds{R}$, we have
\begin{equation}\label{Qb}
- \tilde{C}\ \leqslant\ -\tilde{C}\, \tan\! \left( 3\, \tilde{C}\, t \right)\ \leqslant\ Q(t,Y_{x_2}(t))\ \leqslant\ \tilde{A}\, \tilde{T}\ +\ \tilde{B},   
\end{equation}
which implies that $Q$ cannot blow-up before  $t = \min \left\{ \tilde{T}, T_{max} \right\}$ and $\|Q\|_{L^\infty} \leqslant \sqrt{\tilde{D}/4} (\tilde{T} + 1)$.
Using \eqref{Ricatti1} and the Young inequality $PQ \leqslant {\textstyle \frac{1}{8}}\, P^2 + 2 Q^2 $ one obtains
\begin{align*}
\frac{\mathrm{d}^\lambda}{\mathrm{d}t}\, P (t,X_{x_1}(t))\ 
&\leqslant\ - {\textstyle \frac{1}{4}}\, P(t,X_{x_1}(t))^2\ +\ {\textstyle \frac{19}{8}}\, Q(t,X_{x_1}(t))^2\ +\ 3\, \tilde{C}^2\\
&\leqslant\ - {\textstyle \frac{1}{4}}\, P(t,X_{x_1}(t))^2\ +\ \tilde{D}\, \left( \tilde{T}\ +\ 1 \right)^2\\
&\leqslant\ - {\textstyle \frac{1}{8}}\, P(t,X_{x_1}(t))^2.
\end{align*}
The last inequality follows from $P(t,X_{x_1}(t))^2 \geqslant 4 \tilde{D} ( \tilde{T} + 1 )^2$, which is true initially and holds because the map $t \mapsto P (t,X_{x_1}(t))$ is decreasing and negative. Then $T_{max} \leqslant \tilde{T} \leqslant T$ and from \eqref{P<} we obtain that 
\[
\pushQED{\qed} 
\inf_{[0,T_{max}[ \times \mathds{R}} P(t,x)\ =\ -\infty,
\quad
\sup_{[0,T_{max}[ \times \mathds{R}}  P(t,x)\ <\ +\infty,
\quad
\sup_{[0,T_{max}[ \times \mathds{R}}  |Q(t,x)|\ <\ +\infty.\qedhere
\popQED
\]

\section*{Acknowledgements}
The author would like to thank the anonymous referee for his remarks which help improving the paper.
Many thanks go to Didier Clamond, Denys Dutykh and St\'ephane Junca for the fruitful discussions.


\end{document}